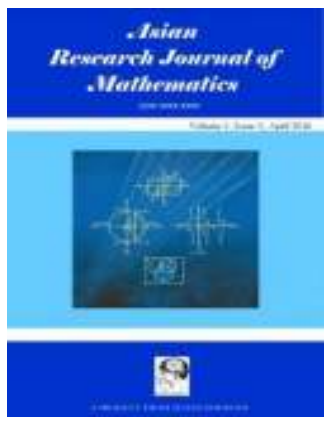

# On Non-existence of Global Weak-predictable-random-field Solutions to a Class of SHEs

**E. M. Omaba[1*], E. Nwaeze[1] and L. O. Omenyi[1]**

*[1]Department of Mathematics, Computer Science, Statistics and Informatics, Faculty of Science, Federal University Ndufu-Alike, Ikwo, P.M.B. 1010, Abakaliki, Ebonyi State, Nigeria.*

***Authors' contributions***

*This work was carried out in collaboration between all authors. All authors read and approved the final manuscript.*

***Article Information***

DOI: 10.9734/ARJOM/2017/33317
*Editor(s):*
(1) Hari Mohan Srivastava, Department of Mathematics and Statistics, University of Victoria, Canada.
*Reviewers:*
(1) Hongwei Zhang, Henan University of Technology, Zhengzhou, China.
(2) Qingying Hu, Henan University of Technology, China.
(3) Pengju Duan, Suzhou University, Anhui, China.
Complete Peer review History: http://www.sciencedomain.org/review-history/19030

***Received: 10th April 2017***
***Accepted: 3rd May 2017***
***Published: 12th May 2017***

*Original Research Article*

## Abstract

The multiplicative non-linearity term is usually assumed to be globally Lipschitz in most results on SPDEs. This work proves that the solutions fail to exist if the non-linearity term grows faster than linear growth. The global non-existence of the solution occurs for some non-linear conditions on $\sigma$. Some precise conditions for existence and uniqueness of the solutions were stated and we have established that the solutions grow in time at most a precise exponential rate at some time interval; and if the solutions satisfy some non-linear conditions then they cease to exist at some finite time $t$. Our result also compares the non-existence of global solutions for both the compensated and non-compensated Poisson noise equations.

---

**Corresponding author: E-mail: mcsylvester.omaba@funai.edu.ng;*

# 1 Introduction

Stochastic PDEs have now received much interest, see [1,2,3] and their references; but the blow up result or rather the non-existence of global result for the discontinuous or jump process has not been given. The study of blow up or global non-existence of solution is as old as Mathematics itself, see [4,5]. The clarity between blow up of a solution and non-existence of global solution of a deterministic PDEs was given by Sugitani, see [6] and the references. Our mild solution here, is a weak-predictable random field solution since it is a class of stochastic heat equation with jump-type process (poisson random noise measure), see [7,8]. We therefore show conditions where the second moment of the solution to the compensated equation and the first moment of the solution to the non-compensated equation (otherwise known as energy solutions) fail to exist and compare the non-existence of global solutions of the two equations with each other. We also consider $\mathcal{L} := -(-\Delta)^{\alpha/2}$, the generator of $\alpha$-stable processes and use some explicit bounds on its corresponding fractional heat kernel to obtain more precise results. We also show that when the solutions satisfy some non-linear growth conditions on $\sigma$, the solutions cease to exist for both compensated and non-compensated noise terms for different conditions on the initial function $u_0(x)$.

Now, consider the following stochastic heat equations driven by a compensated and non-compensated Poisson noises (discontinuous processes).

$$[\frac{\partial u}{\partial t}(t,x) - \mathcal{L}u(t,x)]\mathrm{d}x\mathrm{d}t = \lambda\int_{\mathbf{R}}\sigma(u(t,x),h)\tilde{N}(\mathrm{d}t,\mathrm{d}x,\mathrm{d}h), \tag{1.1}$$

$$[\frac{\partial u}{\partial t}(t,x) - \mathcal{L}u(t,x)]\mathrm{d}x\mathrm{d}t = \lambda\int_{\mathbf{R}^d}\sigma(u(t,x),h)N(\mathrm{d}t,\mathrm{d}x,\mathrm{d}h), \tag{1.2}$$

with initial condition $u(0,x) = u_0(x)$. Here and throughout, $u_0 : \mathbf{R} \to \mathbf{R}_+$ is a non-random function, and $\mathcal{L}$ is the $L^2$- generator of a Lévy process; $\tilde{N}$ is a compensated Poisson random measure and $\lambda$ the noise level.

We follow the method of [9] in making sense of our integral solutions.

**Definition 1.1** *We say that a process $\{u(t,x)\}_{x\in\mathbf{R},t>0}$ is a mild solution of (1.1) if a.s (almost surely), the following is satisfied*

$$\begin{aligned} u(t,x) &= \int_{\mathbf{R}} p(t,x,y)u_0(y)\mathrm{d}y \\ &\quad + \lambda\int_0^t\int_{\mathbf{R}}\int_{\mathbf{R}} p(t-s,x,y)\sigma(u(s,y),h)\tilde{N}(\mathrm{d}h,\mathrm{d}y,\mathrm{d}s), \end{aligned} \tag{1.3}$$

where $p(t,.,.)$ is the heat kernel. If in addition to the above, $\{u(t,x)\}_{x\in\mathbf{R},t>0}$ satisfies the following condition

$$\sup_{0\le t\le T}\sup_{x\in\mathbf{R}} \mathrm{E}\,|u(t,x)|^2 < \infty, \tag{1.4}$$

for all $T > 0$, then we say that $\{u(t,x)\}_{x\in\mathbf{R},t>0}$ is a random field solution to (1.1).

Define

$$\Upsilon(\beta) := \frac{1}{2\pi}\int_{\mathbf{R}} \frac{\mathrm{d}\xi}{\beta + 2\mathrm{Re}\Psi(\xi)} \qquad \text{for all } \beta > 0, \tag{1.5}$$

where $\Psi$ is the characteristic exponent for the Lévy process. A result of Dalang [10] shows that equation (1.1) has a unique solution with the requirement that $\Upsilon(\beta) < \infty$ for all $\beta > 0$ which forces $d = 1$ since $\tilde{N}$ is a martingale valued Poisson measure. Fix some $x_0 \in \mathbf{R}$ and define the *upper* $p$ *th-moment Liapunov exponent* $\bar{\gamma}(p)$ of $u$ at $x_0$ as

$$\bar{\gamma}(p) := \limsup_{t\to\infty} \frac{1}{t} \ln \mathrm{E}\left[ |u(t, x_0)|^p \right] \quad \text{for all } \mathrm{p} \in (0, \infty). \tag{1.6}$$

**Definition 1.2 ( random field solution)**

We seek a mild solution to equation (1.2) of the form.

$$u(t,x) = \int_{\mathbf{R}^d} p(t,x,y)u_0(y)\mathrm{d}y$$
$$+ \lambda \int_0^t \int_{\mathbf{R}^d} \int_{\mathbf{R}^d} p(t-s,x,y)\sigma(u(s,y),h)N(\mathrm{d}h,\mathrm{d}y,\mathrm{d}s), \tag{1.7}$$

with $p(t,.,.)$ the heat kernel. We impose the following integrability condition on the solution:

$$\sup_{t>0} \sup_{x\in\mathbf{R}^d} \mathrm{E}\,|u(t,x)| < \infty.$$

We now make sense of the discontinuous integrals by stating the existence theorem (Ikeda and Watanabe [11]) characterised by a stationary Point process.

**Theorem 1.3** *Given a* $\sigma$*-finite measure* $n(\mathrm{d}x)$ *on* $(X, \mathcal{B}(X))$, *then there exists a stationary Poisson point process* $p$ *if the random measure* $N_p(\mathrm{d}t,\mathrm{d}x)$ *is of the form* $\mathrm{E}[N_p(\mathrm{d}t,\mathrm{d}x)] = n_p(\mathrm{d}t,\mathrm{d}x) = \mathrm{d}t n(\mathrm{d}x)$.

*Proof.* Ikeda and Watanabe [11].

Now applying the above theorem with $X := \mathbf{R}^d \times \mathbf{R}^d$ and $\mathcal{B}(X) := \mathcal{B}(\mathbf{R}^d) \otimes \mathcal{B}(\mathbf{R}^d)$. We will take $n(\mathrm{d}x,\mathrm{d}h) := \mathrm{d}x\nu(\mathrm{d}h)$. One set of the vectors will play the role of position while the other will play the role of "jumps". By the above theorem, we have a Poisson point process $p(s) \in \mathbf{R}^d \times \mathbf{R}^d$. The Poisson random measure is thus given by the following

$$N_p((0,t], A\times B) := \#\{s \le t; s \in D_p; p(s) \in A\times B\},$$

where $D_p$ is a countable subset of $[0,\infty]$.

Let $(\Omega, \mathcal{F}, \{\mathcal{F}_t\}_{t\geq 0}, \mathrm{P})$ be a complete filtered Probability space and $(\mathbf{R}^d, \mathcal{B}(\mathbf{R}^d))$ be a measurable space. Let $\mathcal{F}_t$ be defined by

$$\mathcal{F}_t := \sigma(N_p([0,t], A\times B, \cdot) : A\times B \in \mathcal{B}(\mathbf{R}^d)\times\mathcal{B}(\mathbf{R}^d)) \vee \mathcal{N},$$

where $t > 0$ and $\mathcal{N}$ denotes the null set of $\mathcal{F}$. We can write the Poisson random measure as

$$N_p((0,t], A\times B) := \sum_{s\in D_p, s\leq t} I_{A\times B}(p_x(s), p_h(s)),$$

where we define $p(s) := (p_x(s), p_h(s))$.

In this case, we have $\mathrm{E}[N_p((0,t], A\times B)] = t\,|A|\,\nu(B)$. We now describe the stochastic integral with respect to this Poisson random measure and define the class of integrand precisely:

**Definition 1.4 (The non-compensated Integral)**

$$H_p^1 := \{f(t,x,h) : f \textit{ is } \{\mathcal{F}_t\} - \textit{predictable and } \int_0^t\int_{\mathbf{R}^d}\int_{\mathbf{R}^d} \mathrm{E}\,|f(s,x,h)|\,\mathrm{d}s\mathrm{d}x\nu(\mathrm{d}h) < \infty\}.$$

The following integral can now be defined for all $f \in H_p^1$

$$\int_0^t\int_{\mathbf{R}^d}\int_{\mathbf{R}^d} f(s,x,h,.)N_p(\mathrm{d}s,\mathrm{d}x,\mathrm{d}h) = \sum_{s\leq t, s\in D_p} f(s, p_x(s), p_h(s))$$

as the a.s sum of the following absolutely convergent sum.

**Definition 1.5 (The compensated Integral)** *Define, similarly, for* $f$ *satisfying the square-integrability condition*

$$H_p^2 = \{f(t,x,h) : f \textit{ is } \{\mathcal{F}_t\} - \textit{predictable and } \int_0^t\int_{\mathbf{R}^d}\int_{\mathbf{R}^d} \mathrm{E}\,|f(s,x,h)|^2\,\mathrm{d}s\mathrm{d}x\nu(\mathrm{d}h) < \infty\}.$$

Then for all $f \in H_p^2$, one defines the integral as follows

$$\int_0^t\int_{\mathbf{R}^d}\int_{\mathbf{R}^d} f(s,x,h)\tilde{N}_p(\mathrm{d}s,\mathrm{d}x,\mathrm{d}h) = \sum_{s\leq t, s\in D_p} f(s, p_x(s), p_h(s)) - \int_0^t\int_{\mathbf{R}^d}\int_{\mathbf{R}^d} f(s,x,h)\mathrm{d}s\mathrm{d}x\nu(\mathrm{d}h)$$

as the a.s sum of the following absolutely convergent sum.

For an example, one can define a Poisson random measure $N = \sum_{i\geq 1}\delta_{(T_i, X_i, Z_i)}$ on $\mathbf{R}_+\times\mathbf{R}^d\times\mathbf{R}^d$ defined on a probability space $(\Omega, \mathcal{F}, P)$ with intensity measure $\mathrm{d}t\mathrm{d}x\nu(\mathrm{d}h)$ where $\nu$ is a Lévy measure on $\mathbf{R}^d$; that is, it satisfies the following

$$\int_{\mathbf{R}^d}(1\wedge h^2)\nu(\mathrm{d}h) < \infty.$$

According to [3], let $(\varepsilon_j)_{j\geq 0}$ be a sequence of positive real numbers such that $\varepsilon_j \to 0$ as $j\to\infty$ and $1=\varepsilon_0 > \varepsilon_1 > \varepsilon_2 > \ldots$. Let

$$\Gamma_j = \{h\in \mathbf{R}^d;\, \varepsilon_j <|h|< \varepsilon_{j-1}\},\ j\geq 1 \ \mathit{and}\ \Gamma_0 = \{h\in \mathbf{R}^d; |h|>1\}.$$

Then for any set $B\in \mathcal{B}(\mathbf{R}_+\times\mathbf{R}^d)$, define

$$\int_{B\times\Gamma_j} hN(\mathrm{d}t,\mathrm{d}x,\mathrm{d}h) = \sum_{(T_i,X_i)\in B} Z_i I_{\{Z_i\in\Gamma_j\}},\ j\geq 0,$$

with the following property: $\mathrm{E}[N(B\times\Gamma_j)] = |B|\,\nu(\Gamma_j)$ and

$$N(B\times\Gamma_j) = \#\{i\geq 1; (T_i,X_i,Z_i)\in B\times\Gamma_j\} < \infty\ a.s.$$

Therefore the above integral is finite since the sum contains finitely many terms.

For the existence and uniqueness of (1.1), we need the following condition on $\sigma$. Essentially this condition says that $\sigma$ is globally Lipschitz in the first variable and bounded by another function in the second variable.

**Condition 1.6** *There exist a positive function* $J$ *and a finite positive constant,* $\mathrm{Lip}_\sigma$ *such that for all* $x,y,h\in \mathbf{R}$*, we have*

$$|\sigma(0,h)|\leq J(h) \quad\text{and}\quad |\sigma(x,h)-\sigma(y,h)|\leq J(h)\mathrm{Lip}_\sigma|x-y|. \tag{1.8}$$

The function $J$ is assumed to satisfy the following integrability condition:

$$\int_{\mathbf{R}} J(h)^2\nu(\mathrm{d}h)\leq \mathrm{K}, \tag{1.9}$$

where $\mathrm{K}$ is some finite positive constant.

For the lower bound result, we will need the following extra condition on $\sigma$.

**Condition 1.7** *There exist a positive function* $\bar{J}$ *and a finite positive constant,* $L_\sigma$ *such that for all* $x,h\in \mathbf{R}$*, we have*

$$|\sigma(x,h)|\geq L_\sigma \bar{J}(h)|x| \tag{1.10}$$

The function $\bar{J}$ is assumed to satisfy the following integrability condition:

$$\kappa\leq \int_{\mathbf{R}} \bar{J}(h)^2\nu(\mathrm{d}h)\leq \mathrm{K}, \tag{1.11}$$

where $\mathrm{K}$ is the constant from (1.9) and $\kappa$ is another positive, finite constant.

For the existence and uniqueness of (1.2), we make the following assumption.

**Condition 1.8** *There exist a positive function* $J$ *and a finite positive constant,* $\mathrm{Lip}_\sigma$ *such that for all* $x,y,h \in \mathbf{R}^d$ *, we have*

$$|\sigma(0,h)| \le J(h) \quad \text{and} \quad |\sigma(x,h)-\sigma(y,h)| \le J(h)\mathrm{Lip}_\sigma |x-y|. \tag{1.12}$$

The function $J$ is assumed to satisfy the following integrability condition.

$$\int_{\mathbf{R}^d} J(h)\nu(\mathrm{d}h) \le \mathrm{K}, \tag{1.13}$$

where $\mathrm{K}$ is some finite positive constant.

We can also give the lower bound estimate on the growth of the first moment.

**Condition 1.9** *There exist a positive function* $\bar{J}$ *and a finite positive constant,* $L_\sigma$ *such that for all* $x,h \in \mathbf{R}^d$ *, we have*

$$|\sigma(x,h)| \ge L_\sigma \bar{J}(h)|x| \tag{1.14}$$

The function $\bar{J}$ is assumed to satisfy the following integrability condition.

$$\kappa \le \int_{\mathbf{R}^d} \bar{J}(h)\nu(\mathrm{d}h) \le K, \tag{1.15}$$

where $\mathrm{K}$ is the constant from (1.13) and $\kappa$ is another positive finite constant.

The outline of the paper is given below. The paper is comprised of five sections. Statements and conditions of main results are given in section two. In section three, some estimates on heat kernels are given, some known results and few new propositions are also given. We also give estimates where the proofs of our theorems lie. The proofs of our main results are given in section four and section five contains a brief conclusion of the research.

# 2 Main Results

We show that if the function $\sigma$ grows faster than linear growth, then the energy of the solutions, that is, the second moment $\mathrm{E}|u(t,x)|^2$ of the solution to (1.1) and the first moment $\mathrm{E}|u(t,x)|$ of the solution to (1.2) cease to exist for all time $t$.

Suppose that instead of (1.10), we have the following condition.

**Condition 2.1** *There exists a constant* $\beta > 1$ *such that*

$$|\sigma(x,h)| \ge L_\sigma \bar{J}(h)|x|^\beta, \tag{2.1}$$

where the constant $L_\sigma$ and the function $\bar{J}$ are the same as in condition 1.7.

We then have the following result.

**Theorem 2.2** *Suppose that* $\mathcal{L} := -(-\Delta)^{\alpha/2}$ *and that condition 2.1 is in force. Then there does not exist any random field solution to (1.1).*

If instead of (1.14), we consider the following condition:

**Condition 2.3** *There exists a constant* $\gamma > 1$ *such that*

$$|\sigma(x,h)| \geq L_\sigma \bar{J}(h)\,|x|^\gamma, \tag{2.2}$$

where the constant $L_\sigma$ and the function $J$ are the same as in condition 1.9.

We then have the following result with the initial condition $u_0 : \mathbf{R}^d \to \mathbf{R}_+$, a positive function on a set of positive measure.

**Theorem 2.4** *Suppose that both conditions 1.8 and 2.3 are in force. Then there are no random field solutions to (1.2) whenever the non-negative initial condition* $u_0$ *is bounded below. Let* $\mathcal{L} := -(-\Delta)^{\alpha/2}$, *then under the same conditions, there are no random field solutions to (1.2) even if we only have* $u_0 \neq 0$.

Here, the initial function $u_0$ is assumed to be a bounded non-negative function such that

$$\int_A u_0(x)\mathrm{d}x > 0, \text{ for some } A \subset \mathbf{R}^d.$$

That is, we define $u_0$ as any measurable function $u_0 : \mathbf{R}^d \to \mathbf{R}_+$ which is positive on a set of positive measure. This assumption implies that the set $A = \{x : u_0(x) > \frac{1}{n}\} \subset \mathbf{R}^d$ has positive measure for all but finite many $n$. Thus by Chebyshev's inequality,

$$\int_{\mathbf{R}^d} u_0(x)\mathrm{d}x \geq \int_{\{x:\, u_0(x) > \frac{1}{n}\}} u_0(x)\mathrm{d}x \geq \frac{1}{n}\mu\{x : u_0(x) > \frac{1}{n}\} > 0,$$

where $\mu$ is a Lebesgue measure.

## 2.1 Some auxiliary results and estimates

We present some properties of $p(t,x)$ that will be used in the proof of our results, see [6].

$$\begin{aligned} p(t,x) &= t^{-d/\alpha} p(1, t^{-1/\alpha}x) \\ p(st,x) &= t^{-d/\alpha} p(s, t^{-1/\alpha}x). \end{aligned} \tag{2.3}$$

From the above relation, $p(t,0) = t^{-d/\alpha} p(1,0)$, is a decreasing function of $t$. The heat kernel $p(t,x)$ is also a decreasing function of $|x|$, that's

$$|x| \geq |y| \ \textit{ implies that } \ p(t,x) \leq p(t,y).$$

This and equation (2.3) imply that for all $t \geq s$,

$$p(t,x) = p(t,|x|) = p(s.\frac{t}{s},|x|) = (\frac{t}{s})^{-d/\alpha} p(s,(\frac{t}{s})^{-1/\alpha}|x|)$$
$$\geq (\frac{s}{t})^{d/\alpha} p(s,|x|) \qquad (since \quad (\frac{t}{s})^{-1/\alpha}|x| \leq |x|)$$
$$= (\frac{s}{t})^{d/\alpha} p(s,x).$$

We now state some propositions and lemma whose proofs can be found in [12].

**Proposition 2.5** *[12]. Let* $p(t,x)$ *be the transition density of a strictly* $\alpha$*-stable process. If* $p(t,0) \leq 1$ *and* $a \geq 2$*, then*

$$p(t,\frac{1}{a}(x-y)) \geq p(t,x)p(t,y) \; \forall \, x,y \in \mathbf{R}^d.$$

**Lemma 2.6** *[13,14,15,16]. Suppose that* $p(t,x)$ *denotes the heat kernel for a strictly stable process of order* $\alpha$. *Then the following estimate holds.*

$$p(t,x,y) \approx t^{-d/\alpha} \wedge \frac{t}{|x-y|^{d+\alpha}} \quad \text{for all} \quad t>0 \quad \text{and} \quad x,y \in \mathbf{R}^d.$$

Now we show that the first term $(P_t u_0)(x)$ of the mild solution to (1.1) and (1.2) grow or decay but only polynomially fast with time. Recall that

$$(P_t u_0)(x) := \int_{\mathbf{R}^d} p(t,x,y) u_0(y) \mathrm{d}y.$$

With the assumption that the initial condition $u_0$ is positive on a set of positive measure, we then have the following.

**Proposition 2.7** *[12]. There exist a* $T>0$ *and a positive constant* $c_1$ *such that for all* $t>T$ *and all* $x \in B(0,t^{1/\alpha})$,

$$(P_t u_0)(x) \geq \frac{c_1}{t^{d/\alpha}}.$$

*Proof.* The proof follows by applying Lemma 2.6.

The next proposition is similar to but more general than Proposition 2.7.

**Proposition 2.8** *[12]. Given the assumption on the initial function* $u_0$. *Then for* $t_0 \geq 1$, $\eta > 0$, *there exists* $c(t_0) > 0$ *such that*

$$\int_{\mathbf{R}^d} p(t+t_0,x,y) u_0(y) \mathrm{d}y \geq c(t_0) p(t+\eta, x).$$

Next we give the following a priori result about the continuity of the second moment of the solution to (1.1).

**Proposition 2.9** *[12]. Suppose that condition 1.6 holds, then for each* $x \in \mathbf{R}$ *, the unique solution to (1.1) is mean square continuous in time. That is for each* $x \in \mathbf{R}$ *, the function* $t \mapsto \mathrm{E}[|u(t,x)|^2]$ *is continuous.*

Here also, we present the time continuity of the first moment of the solution (1.2).

**Proposition 2.10** *[12]. Suppose that condition 1.8 holds, then for each* $x \in \mathbf{R}^d$ *, the unique solution to (1.2) is mean continuous in time. That is for each* $x \in \mathbf{R}^d$ *, the function* $t \mapsto \mathrm{E}[|u(t,x)|]$ *is continuous.*

We now give the following proposition which establishes the fact that under the local Lipschitz continuity as stated in condition 2.1, there exists a unique solution up to a fixed time $T$ .

**Proposition 2.11** *[12]. Suppose that condition 2.1 holds. Then there exists a* $T > 0$ *such that (1.1) has a unique random field solution up to time* $T$ *.*

*Proof.* We begin by defining

$$\sigma_N(x,h) = \begin{cases} \sigma(x,h) & \text{if} \quad x \le N \\ \sigma(N,h) & \text{if} \quad x > N. \end{cases}$$

$\sigma_N(x,h)$ therefore satisfies (1.8) but with a different constant. The result follows by the proof of existence and uniqueness of the solution and its time continuity above.

We will also need the following proposition which establishes the fact that under the local Lipschitz continuity as stated in condition 2.2, there exists a unique solution up to a fixed time $T$ .

**Proposition 2.12** *[12]. Suppose that condition 2.2 holds. Then there exists a* $T > 0$ *such that (1.2) has a unique random field solution up to time* $T$ *.*

Here follows Jensen's inequality which shall be used in the proof of the blow-up result.

**Lemma 2.13** *[12] Given that* $p(\mathrm{d}x)$ *is a probability measure on* $\mathbf{R}^d$ *and suppose that the function* $u$ *is non-negative. Then for all convex function* $f$ *the following holds,*

$$\int_{\mathbf{R}^d} f(u(x))\,p(\mathrm{d}x) \ge f\Big(\int_{\mathbf{R}^d} u(x)\,p(\mathrm{d}x)\Big).$$

# 3 Proof of Main Results

In this section, we give the proofs of our main results.

*Proof of Theorem 2.2.* We use Ito's isometry to write

$$\mathrm{E}|u(t,x)|^2 = |(P_t u_0)(x)|^2 + \lambda^2 \int_0^t \int_{\mathbf{R}} \int_{\mathbf{R}} |p(t-s,x-y)|^2\, \mathrm{E}|\sigma(u(s,y),h)|^2\, \nu(\mathrm{d}h)\mathrm{d}y\mathrm{d}s.$$

We use the assumption that $u_0(x) > c_1$ for some positive constant $c_1$ and condition 2.1 to come up with

$$\begin{aligned} \mathrm{E}|u(t,x)|^2 &\geq c_1^2 + \kappa\lambda^2 L_\sigma^2 \int_0^t \int_{\mathbf{R}} p^2(t-s, x-y)\mathrm{E}|u(s,y)|^{2\beta}\, \mathrm{d}y\mathrm{d}s \\ &\geq c_1^2 + \kappa\lambda^2 L_\sigma^2 \int_0^t (\inf_{y\in\mathbf{R}} \mathrm{E}|u(s,y)|^2)^\beta p(2(t-s),0)\mathrm{d}s. \end{aligned}$$

Upon setting

$$F(t) = \inf_{x\in\mathbf{R}} \mathrm{E}|u(t,x)|^2,$$

The above inequality reduces to

$$F(t) \geq c_1^2 + \kappa\lambda^2 L_\sigma^2 \int_0^t p(2(t-s),0)F^\beta(s)\mathrm{d}s.$$

We now use Lemma 2.6 to find lower bounds on the heat kernel appearing in the above display.

$$\begin{aligned} F(t) &\geq c_1^2 + \kappa\lambda^2 L_\sigma^2 c_2 \int_0^t (t-s)^{-1/\alpha} F^\beta(s)\mathrm{d}s \\ &\geq c_1^2 + \kappa\lambda^2 L_\sigma^2 c_2 \int_0^t t^{-1/\alpha} F^\beta(s)\mathrm{d}s. \end{aligned}$$

Hence, multiplying through by $t^{1/\alpha}$

$$\begin{aligned} F(t)t^{\frac{1}{\alpha}} &\geq c_1^2 t^{1/\alpha} + \kappa\lambda^2 L_\sigma^2 c_2 \int_0^t F^\beta(s)\mathrm{d}s \\ &= c_1^2 t^{1/\alpha} + \kappa\lambda^2 L_\sigma^2 c_2 \int_0^t \frac{(s^{1/\alpha}F(s))^\beta}{s^{\beta/\alpha}}\mathrm{d}s. \end{aligned}$$

Let $Y(t) = F(t)t^{1/\alpha}$, then for all $t \geq 0$

$$\begin{aligned} Y(t) &\geq c_1^2 t^{1/\alpha} + \kappa\lambda^2 L_\sigma^2 c_2 \int_0^t \frac{Y^\beta(s)}{s^{\beta/\alpha}}\mathrm{d}s \\ &\geq \kappa\lambda^2 L_\sigma^2 c_2 \int_0^t \frac{Y^\beta(s)}{s^{\beta/\alpha}}\mathrm{d}s, \end{aligned}$$

and solving the ODE: $\dot{Y}(t) = \kappa\lambda^2 L_\sigma^2 c_2 \dfrac{Y^\beta(t)}{t^{\beta/\alpha}}$ with the initial condition $Y(0) = 0$, that is, $\dfrac{\mathrm{d}Y(t)}{Y^\beta(t)} = \kappa\lambda^2 L_\sigma^2 c_2 t^{-\beta/\alpha}\mathrm{d}t$ with $Y(0) = 0$ gives

$$Y(t) = \{\frac{\kappa\lambda^2 L_\sigma^2 c_2 \alpha}{\alpha-\beta} t^{(\alpha-\beta)/\alpha}\}^{1/1-\beta}$$

and $Y(t)$ ceases to exist in finite time for all $\beta$ such that $1 < \beta < \alpha$

We now prove the finite time blow up for the non-compensated noise equation.

*Proof of Theorem 2.4.* We begin with the first part of the theorem. Now write

$$\mathrm{E}|u(t,x)| = |(P_t u_0)(x)| + \lambda \int_0^t \int_{\mathbf{R}^d} \int_{\mathbf{R}^d} |p(t-s,x-y)|\,\mathrm{E}|\sigma(u(s,y),h)|\,\nu(\mathrm{d}h)\mathrm{d}y\mathrm{d}s.$$

We use the assumption that $u_0(x) > c_1$ for some positive constant $c_1$ and condition 2.3 to come up with

$$\begin{aligned}\mathrm{E}|u(t,x)| &\geq c_1 + \kappa\lambda L_\sigma \int_0^t \int_{\mathbf{R}^d} p(t-s,x-y)\mathrm{E}|u(s,y)|^\gamma \,\mathrm{d}y\mathrm{d}s \\ &\geq c_1 + \kappa\lambda L_\sigma \int_0^t (\inf_{y\in\mathbf{R}^d} \mathrm{E}|u(s,y)|)^\gamma \mathrm{d}s.\end{aligned}$$

Upon setting

$$F(t) = \inf_{x\in\mathbf{R}^d} \mathrm{E}|u(t,x)|,$$

the above inequality reduces to

$$F(t) \geq c_1 + \kappa\lambda L_\sigma \int_0^t F^\gamma(s)\mathrm{d}s.$$

We now solve the following differential inequality: $\dot{F}(t) \geq \kappa\lambda L_\sigma F^\gamma(t)$ with $F(0) \geq c_1$ or simply the ordinary differential equation: $\dot{F}(t) = \kappa\lambda L_\sigma F^\gamma(t)$ with the initial condition $F(0) = c_1$, where the solution is given by

$$F(t) = \{(1-\gamma)\kappa\lambda L_\sigma t + c_1^{(1-\gamma)}\}^{1/1-\gamma} = (1-\gamma)^{1/(1-\gamma)}\{\kappa\lambda L_\sigma t + \frac{c_1^{(1-\gamma)}}{1-\gamma}\}^{1/1-\gamma}$$

which fails to converge in a finite time for all $\gamma > 1$ .

Next, we consider when the initial function $u_0$ is positive but not necessarily bounded below and $\mathcal{L} = -(-\Delta)^{\alpha/2}$. Let $t_0 > 0$ be fixed, then

$$\begin{aligned}u(t+t_0,x) &= \int_{\mathbf{R}^d} p(t+t_0,x-y)u(0,y)\mathrm{d}y \\ &\quad + \lambda \int_0^{t+t_0} \int_{\mathbf{R}^d} \int_{\mathbf{R}^d} p(t+t_0-s,x-y)\sigma(u(s,y),h)N(\mathrm{d}h,\mathrm{d}y,\mathrm{d}s).\end{aligned}$$

Taking first moment of both sides, then by Proposition 2.8 we have

$$
\begin{aligned}
\mathrm{E}|u(t+t_0,x)| &\geq \int_{\mathbf{R}^d} p(t+t_0,x-y)u(0,y)\mathrm{d}y \\
&+\lambda\int_0^{t+t_0}\int_{\mathbf{R}^d}\int_{\mathbf{R}^d} p(t+t_0-s,x-y)\mathrm{E}|\sigma(u(s,y),h)|\nu(\mathrm{d}h)\mathrm{d}y\mathrm{d}s \\
&\geq c(t_0)p(t+\eta,x) \\
&+\lambda L_\sigma \kappa\int_{t_0}^{t+t_0}\int_{\mathbf{R}^d} p(t+t_0-s,x-y)\mathrm{E}|u(s,y)|^\gamma \mathrm{d}y\mathrm{d}s.
\end{aligned}
$$

We perform a change of variable, by letting $\tau = s-t_0$, $\mathrm{d}\tau = \mathrm{d}s$ and therefore

$$
\begin{aligned}
\mathrm{E}|u(t+t_0,x)| &\geq c(t_0)p(t+\eta,x) \\
&+\lambda L_\sigma \kappa\int_0^t\int_{\mathbf{R}^d} p(t-\tau,x-y)\mathrm{E}|u(\tau+t_0,y)|^\gamma \mathrm{d}y\mathrm{d}\tau.
\end{aligned}
$$

Let $v(t,x) := u(t+t_0,x)$; then to show that $\mathrm{E}|u(t+t_0,x)|$ fails to exist in some finite time suffices to show the same for $\mathrm{E}|v(t,x)|$ and therefore

$$
\mathrm{E}|v(t,x)| \geq c(t_0)p(t+\eta,x)+\lambda L_\sigma\kappa\int_0^t\int_{\mathbf{R}^d} p(t-\tau,x-y)(\mathrm{E}|v(\tau,y)|)^\gamma \mathrm{d}y\mathrm{d}\tau.
$$

Multiplying through by $p(t,x)$, and integrating in $\mathrm{d}x$, we obtain

$$
\begin{aligned}
\int_{\mathbf{R}^d}\mathrm{E}|v(t,x)|\,p(t,x)\mathrm{d}x &\geq c(t_0)\int_{\mathbf{R}^d} p(t+\eta,x)p(t,x)\mathrm{d}x \\
&+\lambda L_\sigma\kappa\int_0^t \mathrm{d}\tau\int_{\mathbf{R}^d}\mathrm{d}y(\mathrm{E}|v(\tau,y)|)^\gamma\int_{\mathbf{R}^d} p(t-\tau,x-y)p(t,x)\mathrm{d}x \\
&= c(t_0)p(2t+\eta,0)+\lambda L_\sigma\kappa\int_0^t \mathrm{d}\tau\int_{\mathbf{R}^d}\mathrm{d}y(\mathrm{E}|v(\tau,y)|)^\gamma p(2t-\tau,y).
\end{aligned}
$$

The last line follows by Kolmogorov property. Moreover, by properties of $p(t,x)$, Jensen's inequality and Lemma 2.13, we get

$$
\begin{aligned}
&\int_{\mathbf{R}^d}\mathrm{E}|v(t,x)|\,p(t,x)\mathrm{d}x \\
&\geq c(t_0)p(1,0)(2t+\eta)^{-\frac{d}{\alpha}}+\kappa L_\sigma\kappa\int_0^t \mathrm{d}\tau\int_{\mathbf{R}^d}\mathrm{d}y\left(\frac{\tau}{2t-\tau}\right)^{\frac{d}{\alpha}} p(\tau,y)(\mathrm{E}|v(\tau,y)|)^\gamma \\
&\geq c_0(2t+\eta)^{-\frac{d}{\alpha}}+\lambda L_\sigma\kappa\int_0^t \mathrm{d}\tau\left(\frac{\tau}{2t-\tau}\right)^{\frac{d}{\alpha}}\left(\int_{\mathbf{R}^d}(\mathrm{E}|v(\tau,y)|\,p(\tau,y)\mathrm{d}y\right)^\gamma.
\end{aligned}
$$

Set $F(t) = \int_{\mathbf{R}^d}\mathrm{E}|v(t,x)|\,p(t,x)\mathrm{d}x$ then

$$
\begin{aligned}
F(t) &\geq c_0(2t+\eta)^{-\frac{d}{\alpha}}+\lambda L_\sigma\kappa\int_0^t\left(\frac{\tau}{2t-\tau}\right)^{\frac{d}{\alpha}}F(\tau)^\gamma \mathrm{d}\tau \\
&\geq c_0(2t+\eta)^{-\frac{d}{\alpha}}+\lambda L_\sigma\kappa\int_0^t\left(\frac{\tau}{2t}\right)^{\frac{d}{\alpha}}F(\tau)^\gamma \mathrm{d}\tau.
\end{aligned}
$$

Multiply through by $t^{\frac{d}{\alpha}}$, therefore

$$F(t)t^{\frac{d}{\alpha}} \geq c_0(\frac{t}{2t+\eta})^{d/\alpha} + \frac{\lambda L_\sigma \kappa}{2^{d/\alpha}}\int_0^t \frac{(\tau^{\frac{d}{\alpha}}F(\tau))^\gamma}{\tau^{\frac{d(\gamma-1)}{\alpha}}}\mathrm{d}\tau.$$

Let $Y(t) = F(t)t^{\frac{d}{\alpha}}$ and assume further that $t \geq \delta$ for all $\delta \geq 0,$ then

$$Y(t) \geq c_0(\frac{\delta}{2\delta+\eta})^{d/\alpha} + \frac{\lambda L_\sigma \kappa}{2^{d/\alpha}}\int_\delta^t \frac{Y(\tau)^\gamma}{\tau^{\frac{d(\gamma-1)}{\alpha}}}\mathrm{d}\tau \geq \frac{\lambda L_\sigma \kappa}{2^{d/\alpha}}\int_\delta^t \frac{Y(\tau)^\gamma}{\tau^{\frac{d(\gamma-1)}{\alpha}}}\mathrm{d}\tau$$

and the result follows by similar argument as the above proof for all $\gamma$ such that $1 < \gamma < 1+\alpha/d$

# 4 Conclusion

In conclusion, we see that both (1.1) and (1.2) fail to have global solutions when the initial condition $u_0$ is bounded below for the operator $\mathcal{L}$ - the $L^2$ generator of a Lévy process; and for the fractional Laplacian $\mathcal{L} = (-\Delta)^{\alpha/2}, \alpha \in (1,2]$ the solution to (1.2) fails to exist globally when the initial function $u_0$ is a positive function on a set of positive measure.

# Acknowledgement

Our appreciation goes to the anonymous referees for thoroughly and carefully reading through the manuscript thereby making good observations that improved the content of the paper.

# Competing Interests

Authors have declared that no competing interests exist.

---